\documentclass{birkjour}
\usepackage{mathptmx}      
\usepackage{amssymb}
 \theoremstyle{definition}
 
 \theoremstyle{remark}

 \numberwithin{equation}{section}
 
\usepackage{orcidlink}

\usepackage{xcolor,soul}
\newcommand{\sbytext}[1]{{ #1}}

\begin{document}

%
%
%
%
%
%
%
%
%

\title[An Adaptive Difference Method]
 {An Adaptive Difference Method for \\ Variable-Order Diffusion Equations}

\author[J. Quintana-Murillo]{Joaqu\'{\i}n Quintana-Murillo}

\address{%
Departamento de Matem\'aticas Universidad de C\'ordoba\\
Campus de Rabanales Edificio Albert Einstein\\
 E-14071  C\'{o}rdoba\\
Spain}
\email{jquintanam01@educarex.es}

\thanks{Financial support from the Grant No.\ PID2020-112936GB-I00/AEI/10.13039/501100011033 and from the Junta de Extremadura (Spain) through Grant Nos.\ IB20079 and GR18079, all of them partially financed by Fondo Europeo de Desarrollo Regional funds.}
\author[S. B. Yuste]{Santos Bravo Yuste \orcidlink{0000-0001-8679-4195} }
\address{Departamento de F\'{\i}sica Universidad de Extremadura\\
Avenida de Elvas, s/n\\
E-06071 Badajoz\\
Spain}
\email{santos@unex.es}

\keywords{Variable-order equations, Fractional derivative, Adaptive difference method, Caputo derivative}


\begin{abstract}
An adaptive finite difference scheme for variable-order fractional-time subdiffusion equations in the Caputo form is studied. The fractional time derivative is discretized by  the L1 procedure but  using nonhomogeneous timesteps.  The size of these timesteps is chosen by  an adaptive algorithm in order to keep the local error bounded around a preset value, a value that can be chosen at will. For some types of problems this adaptive method is much faster than the corresponding usual method with fixed timesteps while keeping the local error of the numerical solution around the preset values. These findings turns out to be similar to those found for constant-order fractional diffusion equations.
\end{abstract}

\maketitle



\section{Introduction}
\label{secIntro}
Fractional calculus has become a useful branch of Mathematics with a wide range of applications in Science and Engineering \cite{Klafter2011,Klages2008,Sokolov2002}. In particular, fractional diffusion equations appears naturally as a useful way of describing some stochastic processes leading to anomalous diffusion \cite{Metzler2000,Sokolov2005,Sokolov2002}. When the diffusion process is normal, the mean square displacement $\langle x^2\rangle$ of the diffusive entity (``particle'' or ``walker'')  is proportional to the time $t$. But in many instances of Nature and social systems, one finds that   $\langle x^2\rangle \propto t^\gamma$ with $\gamma\neq 1$.  In this case the diffusion is anomalous. There is subdiffusion when $\gamma<1$ and superdiffusion when $\gamma>1$. The value of the anomalous diffusion exponent $\gamma$ depends on the type of particle and medium in which the particle moves. For example, one finds anomalous diffusion when a diffusing particle, a random walker, travels on disordered media \cite{Havlin1987,Klages2008}.  Also, one finds anomalous diffusion when the time  $\tau$ between steps of the random walker (even in a normal, no fractal, Euclidean medium) is a random variable drawn from a heavy tail distribution $\varphi(\tau)\sim \tau^{-1-\gamma}$ with $0<\gamma<1$. This class of diffusive systems is described by the so-called Continuos Time Random Walk (CTRW) model which leads naturally to the fact that the probability density $u(\vec r,t)$ of finding a particle at $\vec r$ at time $t$ follows some class of  fractional diffusion equation.
For example, for a large class of random walks in one dimension with no reactions or source terms, the fractional partial differential equation (FPDE)  for $u(x,t)$ is
\begin{equation}
\label{difuEqRL}
\frac{\partial }{\partial t}u(x,t)=\;_{0}D_t^{1-\gamma} \mathcal{L}_\text{FP} \,  u(x,t)
\end{equation}
where
$\,_{0}D_t^{1-\gamma} $ is the Riemann-Liouville derivative defined by
\begin{equation}
\label{RLderivative}
\;_{0}D_t^{1-\gamma} f(t)=
\frac{1}{\Gamma(\gamma)} \frac{d}{dt}\int_0^t d\tau  \frac{f(\tau)}{(t-\tau)^{1-\gamma}}
\end{equation}
and $\mathcal{L}_\text{FP}$ is the Fokker-Planck operator \cite{LeVot2018,Metzler1999b,Metzler2000}
\begin{equation}
\label{FPoperator}
\mathcal{L}_\text{FP} = K_\gamma \frac{\partial^2  }{\partial x^2} -\frac{\partial}{\partial x}v_\gamma(x,t).
\end{equation}
Here $K (\partial^2/\partial x^2)$ is the diffusion term,  $ (\partial/\partial x)v_\gamma(x,t)$ is the advection term,  $K $ is the anomalous diffusion coefficient, and  $v_\gamma(x,t)$ is related to the external force applied to the random walker.  When $v(x,t)=0$ one has  $\langle x^2(t)\rangle \sim K  t^\gamma$ for large $t$, being $\langle x^2(t)\rangle$ the mean square displacement of the walker.
Strictly, the right fractional operator in Eq.~\eqref{difuEqRL} is not the Riemann-Liouville derivative but the Gr{\"{u}}nwald-Letnikov derivative. However, both operators are equivalent for well-behaved solutions $f(t)$ where $\lim_{t\to 0}\int_0^t
d\tau\,(t-\tau)^{\gamma-1}f(\tau) =0$ \cite{Podlubny1999}. This property is usually assumed.
An alternative way of writing Eq.~\eqref{difuEqRL} is
\begin{equation}
\label{difuEqCaputo}
  \frac{\partial^{\gamma}}{\partial t^{\gamma}}u(x,t)=
\mathcal{L}_\text{FP} u(x,t)
\end{equation}
where
\begin{equation}
\label{CaputoDerivative}
  \frac{\partial^{\gamma}}{\partial t^\gamma}f(t)=
\frac{1}{\Gamma(1-\gamma)} \int_0^t d\tau  \frac{1}{(t-\tau)^{\gamma}}
\frac{df(\tau)}{d\tau}
\end{equation}
is the Caputo derivative. When $\gamma$ is constant, Eqs.~\eqref{difuEqRL} and \eqref{difuEqCaputo} (in the language of Ref.\cite{Sokolov2005}, the ``modified'' and ``normal'' form of the FPDE, respectively)  are equivalent.
\sbytext{Since the value of the anomalous diffusion exponent $\gamma$ depends on the type of medium in which diffusion occurs, $\gamma$ becomes a function of space and/or time in systems where the medium is spatially and/or temporally heterogeneous} \cite{Chechkin2005,Fedotov2019,Fedotov2021,Sun2019}.
In this case the FPDE equation is called variable order FPDE (VOFPDE).  The equation is usually called fractional diffusion (or subdiffusion) equation when  $0<\gamma(x,t)<1$. This is the equation on which we are going to focus, although our  procedure is also valid for the so-called diffusion-wave equation in which  $1<\gamma(x,t)<2$ \cite{Mainardi1996,Schneider1989}.

As for standard partial differential equations, it is no possible to find analytical solutions for the FPDEs in many cases (especially for VOFPDEs)  and one has then  to resort to numerical methods.  An important class of these methods are finite difference methods \cite{Li2018,Li2015}. Some comments about these methods when applied to FPDE are in order.
First, let $N$ be the number of timesteps of size $\tau_i$ required by the finite difference method to provide the solution $u(x,t)$  at time $t$, i.e. $t=\sum_{i=1}^N\tau_i$.  The computation time required by these methods is known to scale as $N^2$, which calls for methods for which $N$ be as smaller as possible.  A way to reach this goal is by using a high order method so that its accuracy is good even for relatively large values of the timestep; see, e.g.,  Refs.~\cite{Ford2017,Hao2017,Li2018,Sun2021,Wang2020,Xing2018,Zeng2015} and references there in.
Second, the singular nature of the integro-differential derivatives around $t=0$ \emph{often} manifests in the rapid change of the solution $u(x,t)$ around this point, specifically, it manifests in the singularity of the first time derivatives of $u(x,t)$ around $t=0$. A way of seeing this is recalling that the solution of Eqs. \eqref{difuEqRL} and  \eqref{difuEqCaputo} can be often written as superposition of subdiffusive modes $u(x,t)=\sum \phi_n(x) E_\gamma(-c_n\, t^\gamma)$  where $\phi_n(x)$ are the eigenfunctions of the Fokker-Planck operator $\mathcal{L}_{FP}$ and $E_\gamma(\cdot)$ \sbytext{is the (one-parameter)}  Mittag-Leffler function \cite{Metzler2000,Podlubny2013}.
For normal diffusion, $\gamma=1$, the Mittag-Leffer functions become exponentials and then one recover the standard expansion of the solution $u(x,t)$ in terms of normal diffusive modes \cite{Carslaw1959}. But, unlike the exponential function, the Mittag-Leffler function has singular derivatives at $t=0$.
This is easy to see recalling that $E_\gamma(z)=\sum_{k=0}^\infty z^k/\Gamma(1+ \gamma k)$, which implies that the $n$-th derivative of $E_\gamma(-c_n\, t^\gamma)$ goes as $t^{\gamma-n}$ for short times (see more details in \cite{Stynes2017}).
This implies that  the solution $u(x,t)$ changes very fast for short times.  This generic behaviour of the solutions of the FPDEs is often overlooked when the numerical methods are tested because often the FPDE equations chosen for this task are carefully build adding (often awkward) extra terms (force terms) so that the equations have simple solutions (usually in the form of a polinomial in $t$).
Unfortunately, in many cases this is the price one has to pay in order to be able to compare numerical solutions with exact solutions.
Finally, note that that $E_\gamma(-\text{cte}\, t^\gamma) \sim t^{-\gamma}$ for large $t$, which implies that \emph{often} the solution $u(x,t)$ of FPDE changes very slowly for large times.

Summing up,  typically the problems we face involve at least two quite different timescales as the solution $u(x,t)$ of the fractional diffusion equations changes quite rapidly for short times and very slowly for long times.
This calls for the adaptation of the size of the timesteps $\tau_i$ to this behaviour: one should choose small timesteps when the solution changes rapidly in order to follow the changes in the solutions and large timesteps when  the solution evolves slowly.  This flexibility in choosing the sizes of the time intervals at our disposal can be exploited in different ways, in particular, it can be used to follow in detail the evolution of the solution in a certain time interval that interest us. This procedure has been recently employed for FPDE of constant order  \cite{Fan2019,Fazio2018,Jannelli2020,Li2016,Podlubny2013,Quintana-Murillo2013,Stynes2017,Yuste2012,Yuste2016,Zhang2014,Zhang2021}. In a particular subclass of these methods the choice of the stepsizes is made taking into account the behaviour of the solution on the fly, that is, taking into account the behaviour of the solution that the numerical method is obtaining in every timestep \cite{Jannelli2020,Yuste2016}.  In this paper we  explore this approach, i.e.,  variable timesteps plus an adaptive algorithm,  for variable-order FPDEs (VOFPDEs) of the form
\begin{equation}\label{FDEeq}
\partial u=F(x,t)
 \end{equation}
with  $\partial\equiv \partial^\gamma/\partial t^\gamma -\mathcal{L}_\text{FP}$  and $\gamma=\gamma(x,t)$.
Therefore, this paper can then be seen as a generalization to variable-order FPDEs of the procedures put forward in Refs.~\cite{Yuste2012}  and \cite{Yuste2016} for constant-order FPDEs.  Some recent numerical methods for solving VOFPDEs   can be found in Refs. \cite{Adel2018,Cao2017,Hajipour2018,Liao2014}; see also the recent reviews by  Sun et al.~\cite{Sun2019} and Patnaik et al. \cite{Patnaik2020} on physical models, numerical methods and applications of VOFDEs.

Finally, what is relevant in this article is how we deal with the fractional operator of Eq.~\eqref{FDEeq}; our treatment of the Fokker-Planck operator is fairly standard. Therefore, for the sake of simplicity, we will not consider the advection term in the Fokker-Planck operator (although its inclusion would be straightforward).
Thus, the equation we will consider in this paper is Eq.~\eqref{FDEeq} with $\mathcal{L}_\text{FP}= K   \partial^2/\partial x^2$:
\begin{equation}\label{eqGraltoSolve}
\frac{\partial^{\gamma(x,t)}}{\partial t^{\gamma(x,t)}}u(x,t)=
  K  \frac{\partial^2}{\partial x^2} u(x,t)+F(x,t).
\end{equation}

This paper is organized as follows. In  Sect.\ \ref{secNumAlgo} we present the numerical scheme and the adaptive algorithm we use.  In  Sect.\  \ref{secTests}  we discuss some of the  properties and qualities of our adaptive method by applying it to several problems.  The paper is closed in Sect.\ \ref{secSumConclu} with some concluding remarks.

\section{The numerical scheme}
\label{secNumAlgo}

The numerical scheme we discuss in this paper has two main components: a difference method able to work with non-uniform timesteps and an algorithm that adapts the size of the timesteps to the behaviour of the solution.  We next present the difference method  and left the discussion of the adaptive algorithm for the final part of this section.

\subsection{The  difference method}
\label{secsubNumDiffMet}
The  difference method we use in this paper is a straightforward generalization to variable order FDEs of the method discussed in \cite{Yuste2012,Zhang2014} for FDEs of constant order. It will be obtained by discretizing the operator  $\partial$ of Eq.~\eqref{FDEeq} by means of the three-point centered formula of the Laplacian and by means of a generalization of the L1 formula for non-uniform meshes of the Caputo derivative.

As is standard for finite difference methods, one looks for the solution $u(x_j,t_m)\equiv u_j^m$ of the continuous integro-differential equation at the nodes $(x_j , t_m)$ of a mesh that cover the space-time region where one wants to find the approximate numerical solution. We start discretizing the Caputo derivative:
\begin{align}
\left.\frac{\partial^\gamma u}{\partial t^\gamma}\right\rvert_{(x_j,t_n)}
&=
\frac{1}{\Gamma(1-\gamma_j^n)} \sum_{m=0}^{n-1} \int_{t_m}^{t_{m+1}} \frac{dt'}{(t-t')^{\gamma_j^n}} \,\frac{\partial u}{\partial t'}
\nonumber\\
&=
\frac{1}{\Gamma(1-\gamma_j^n)} \sum_{m=0}^{n-1} \,\frac{u_j^{m+1}-u_j^m}{\tau_{m+1}} \int_{t_m}^{t_{m+1}} \frac{dt'}{(t-t')^{\gamma_j^n}}+R_j(t_n)  \nonumber\\
&=
\frac{1}{\Gamma(2-\gamma_j^n)} \sum_{m=0}^{n-1} T_j^{m,n}\,\frac{u_j^{m+1}-u_j^m}{\tau_{m+1}} +R_j(t_n)
\label{discreCaputo}
\end{align}
where $t_0=0$, $n\ge 1$, $\tau_m\equiv t_{m}-t_{m-1}$,  $\gamma_j^n\equiv \gamma(x_j,t_n)$,
\begin{equation}
\label{x}
T_j^{m,n}=\frac{(t_n-t_m)^{1-\gamma_j^n}-(t_n-t_{m+1})^{1-\gamma_j^n}}{\tau_{m+1}},\qquad m\le n-1,
\end{equation}
$T_j^{0,1}=(t_1-t_0)^{-\gamma_j^1}$, and $R_j(t_n)$ is the temporal truncation error. This is generalization of the L1 formula \cite{Gao2014,Oldham1974} for non-uniform meshes of the Caputo derivative. It is easy to see \cite{Yuste2012} that $R_j(t_n)$ is bounded by a term of order $t_n^{1-\gamma_j^n} \tau_\text{max}$ with $\tau_\text{max}=\max_{1\le m\le n} \tau_{m+1}$. A improved bound of $R_j(t_n)$  when $u(\cdot,t)\in \mathcal{C}^2[0,t_n]$ is $[\tau_n^2/(2-2\gamma)+\tau_\text{max}^2/8)]\tau_n^{-\gamma} \max_{0\le t\le t_n} \partial^2 u(\cdot,t)/\partial t^2$ \cite{Zhang2014}.

Regarding the spatial discretization, we use a mesh of fixed size, $x_{j+1}-x_j=\Delta x$ for all $j$, and we discretize the Laplacian operator $\partial^2/\partial x^2$   by means of the three-point centered formula:
\begin{equation}
\label{Lapladiscre2}
\left. \frac{\partial^2 u(x,t)}{\partial x^2}\right\rvert_{(x_j,t_n)}=\frac{u_{j+1}^n-2u_j^n+u_{j-1}^n}{(\Delta x)^2}+R^n(x_j)
\end{equation}
with $R^n(x_j)=O(\Delta x)^2$.

Introducing \eqref{Lapladiscre2} and \eqref{discreCaputo} in the equation  that we want to solve, Eq.~\eqref{eqGraltoSolve}, neglecting the truncation errors $R^n(x_j) $ and $R_j(t_n)$, and reordering the terms, we obtain the finite difference scheme we were looking for:
   \begin{equation}
   \label{metImp1}
  -S_j^n \, U^n_{j+1} +(1+2S_j^n) U^n_{j}-S_j^n \, U^n_{j-1} =\mathcal{M}\left[U_j^{(n)}\right]+ \tilde F_j^n
\end{equation}
where
 \begin{align}
\label{Sjn}
S_j^n&= \Gamma(2-\gamma_j^n)  K \frac{ (t_n-t_{n-1})^{\gamma_j^n}}{(\Delta x)^2}.
\\
\label{Moperator}
\mathcal{M}\left[U_j^{(n)}\right]&= U^{n-1}_{j} - \sum_{m=0}^{n-2} \tilde T_j^{m,n} \left[ U^{m+1}_{j} -U^{m}_{j}\right],
\end{align}
$\tilde T_j^{m,n}=(t_n-t_{n-1})^{\gamma_j^n} T_j^{m,n}$, $\tilde F_j^n=(t_n-t_{n-1})^{\gamma_j^n}  F(x_j,t_n)$, and $U^n_{j}$ is the numerical estimate of the exact value $u_j^n$.  This implicit finite difference scheme is in the form of a tridiagonal system and can therefore be solved easily and efficiently.
It is also unconditionally stable and convergent. This can be proved by means of a straightforward extension of the arguments used in \cite{Yuste2012} (and also in \cite{Zhang2014}) for the case $\gamma(x,t)=\text{cte.}$ In fact, if one compares the present scheme \eqref{metImp1} with that of Ref.~\cite{Yuste2012}, one sees that the only difference is that the quantities $S_j^n$ and $T_j^{m,n}$ depend here on $j$. But, for our purposes, this fact is not relevant and so the  procedures employed in Ref.~\cite{Yuste2012} can be applied directly to the present variable-order scheme  to show that it is convergent and unconditionally stable. We will not repeat these arguments here.

\subsection{The adaptive algorithm}
\label{secsubAdapAlgo}

In this paper, the dynamical way (the adaptive algorithm) we use to choose the size of the timesteps is a step-doubling algorithm \cite{Press2007,Quintana-Murillo2013} that we called the ``trial \& error'' method in \cite{Yuste2016}. Assuming the the values of $U_k^{(m)}$ with $0<m<n-1$ have been already evaluated, in this method (see more details in \cite{Quintana-Murillo2013,Yuste2016}):
\begin{itemize}
\item The solutions at the time $t_n$ of the $n$-th timestep are evaluated twice, first one gets the solution $U_k^{(n)}$ employing a full step $\Delta_n=t_n-t_{n-1}$ and then the solution $\hat{U}_k^{(n)}$ using two half steps of size $\Delta_n/2$.
\item The ``difference'' $\mathcal{E}^{(n)}=\underset{\text{all} k}{\text{max}} \left\lvert U_k^{(m)}-\hat{U}_k^{(n)}\right\rvert$  between these two numerical estimates of the solution is used as an indicator of the local error of the numerical method.
\item The final size of $\mathcal{E}^{(n)}$ is kept around a prefixed tolerance $\tau$ by adjusting the value of $\Delta_n$:   if $\mathcal{E}^{(n)}$ is larger that the tolerance, we halve the size of $\Delta_n$ until $\mathcal{E}^{(n)}<\tau$, but if $\mathcal{E}^{(n)}$ is smaller that the tolerance, we double the size of $\Delta_n$ until $\mathcal{E}^{(n)}>\tau$, and we take as the final $\Delta_n$ the last one for which $\mathcal{E}^{(n)}$ was smaller than the tolerance.
\end{itemize}

Other adaptive algorithms that work well for non-fractional (ordinary) problems could be extended to our FPDE. For example, predictive methods, based on the (often assumed) dependence of local error on the size of timesteps, have been extended to FPDE \cite{Yuste2016}.   Recently, Jannelli \cite{Jannelli2020} has proposed an adaptive  procedure for FPDE that is a generalization of an adaptive procedure developed for stiff problems described by ordinary differential equations.  We might also be interested in an accurate evaluation of the solution within a given time interval (or time regime) where some quantity (say the variable order exponent) exhibits a certain property; in this case, the adaptive algorithm should adapt the size of the timesteps to carefully explore this time interval. Finally, one might be more interested in accurately evaluating some quantity (a conserved quantity, for example) than in evaluating the solution; the adaptive algorithm should then be designed accordingly.

\section{Results}
\label{secTests}

In this section we are going to discuss some of the properties and  advantages (and disadvantages) of our adaptive method applying it to four different problems with VOFDEs. In the first example, the fractional derivative depends only on time and has a simple exact solution. This case will serve as testbed to study the precision and efficiency (required CPU time) of the method. In the second and third example, $\gamma$ only depends on the position $x$.  The second example illustrates the case of a sudden change in the solution for small $t$, while in the third example we study a problem with a solution that goes towards a nonzero stationary solution.   In the fourth and final example, $\gamma$ is a periodic function similar to a pulsating function. This example  shows how the adaptive method adapts by periodically changing the size of the time intervals.

\subsection{Case 1}

The first problem we consider has a very simple and well-behaved solution, namely,
\begin{equation}
\label{solCase1}
u(x,t)= \left(2-e^{-t}\right)\sin x.
\end{equation}
Taking into account that \cite{Garrappa2019}
\begin{equation}
\label{x}
  \frac{\partial^{\gamma}}{\partial t^\gamma} e^{\lambda t}= \lambda\, t^{1-\gamma}\, E_{1,2-\gamma}(\lambda t)
\end{equation}
for $0<\gamma\le 1$, one finds that Eq.~\eqref{solCase1} is solution of the FPDE
\begin{subequations}
\label{eqsCaso1}
\begin{align}
        \frac{\partial^{\gamma}u(x,t)}{\partial t^\gamma} &=\frac{\partial^2 u(x,t)}{\partial x^2}+F(x,t), \label{eqCaso1} \\
        u(x,0)&=\sin x,\qquad 0\leq x \leq \pi, \label{eqCaso1CI} \\
        u(0,t)&=u(\pi,t)=0, \label{eqCaso1CC} \\
        F(x,t)&=\left[2-e^{-t}+t^{1-\gamma}\, E_{1,2-\gamma}(-t)\right]\,\sin x,  \label{eqCaso1TF}
\end{align}
\end{subequations}
where $E_{1,2-\gamma}(\cdot)$ is the two-parameter  Mittag-Leffler function \cite{Podlubny1999}.  In this example we consider the VOFPDE in which $\gamma\equiv\gamma(t)=\left(1+ e^{-t}\right)/2$, that is, a VOFPDE where the equation goes from a standard, non-fractional,  PDE with $\gamma=1$ for $t=0$ to a strongly anomalous diffusive equation with $\gamma=1/2$ for $t\to\infty$.

\begin{figure}
\begin{center}
\includegraphics[width=0.75\textwidth]{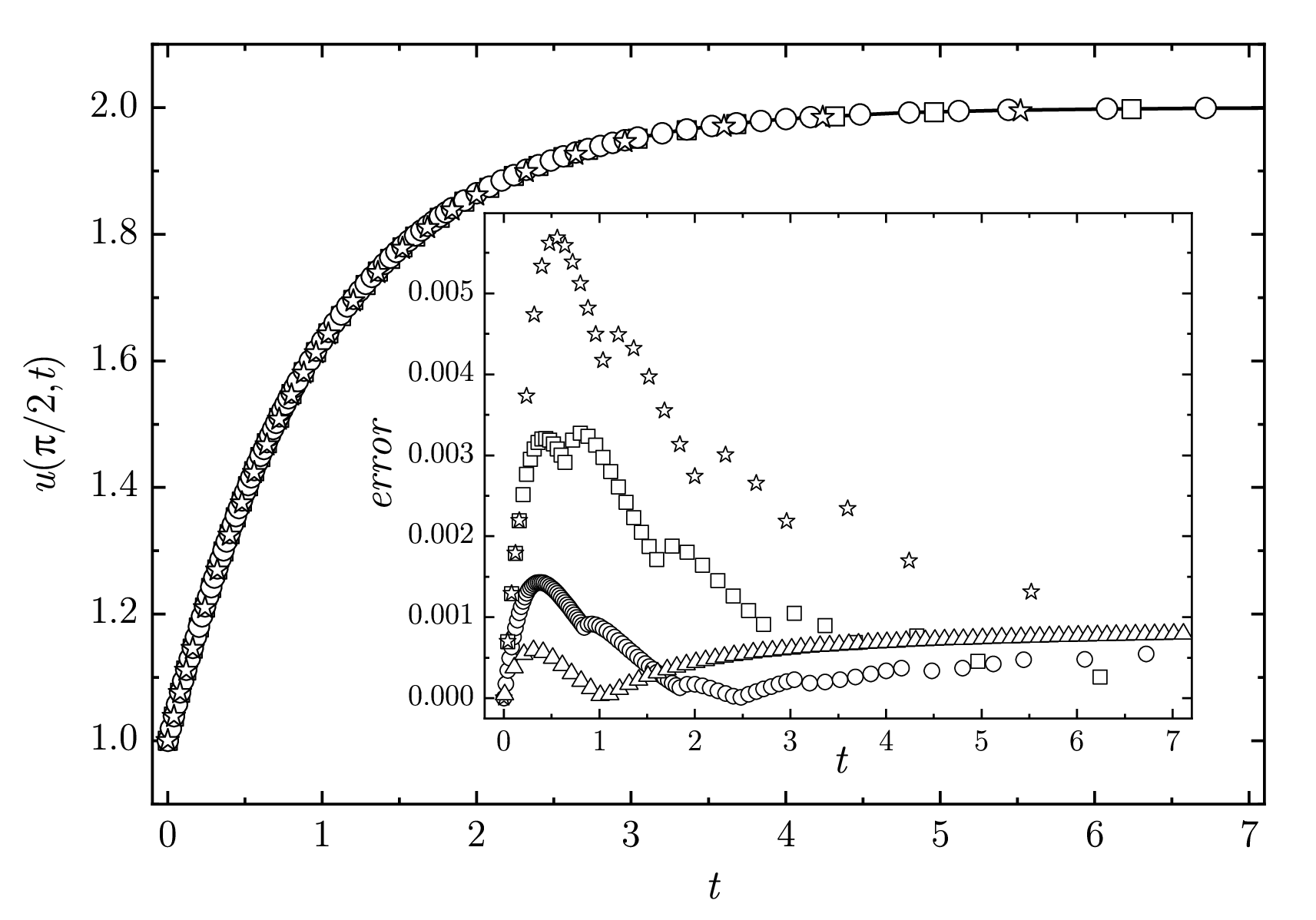}    
\caption{$u(\pi/2, t)$ vs. $t$ for the problem \eqref{eqsCaso1} of Case 1. The symbols are the numerical solutions provided by the adaptive method with tolerances $\tau=10^{-3}$ (stars), $\tau=5\times 10^{-3}$ (squares) and $\tau=10^{-4}$ (circles). In all cases $\Delta x=\pi/40$.
The line is the exact analytical solution.
The numerical errors $|U_j^n-u(x_j=1/2,t_n)|$ are plotted in the inset.  As reference, we have also included the numerical errors when fixed timesteps of size 0.01 are used (triangles). For the sake of readability, we have only plotted one of every 10 points in this case.
}
\label{fig:Ucase2}
\end{center}
\end{figure}

In Fig.~\ref{fig:Ucase2} we compare the  analytical solution (solid line) at the middle point, $u(\pi/2, t)$,   with the numerical solution for three different values of the tolerance $\tau$.
We see that the agreement of the numerical solutions with the exact one is excellent.  We also note that, as expected,  the size of the timesteps changes according to the behavior of the solution: for short times,  the timesteps are small because the solution changes rapidly, which leads to the noticeable cluttering of points in this time range; on the other hand, for large times  the timesteps are very large because the solution hardly changes.   In fact, we know that $u(x=1/2,t)$ goes asymptotically towards the constant value $2$ for long times. This implies that the numerical estimate of the solution will be very accurate even for very large values of the timesteps.  This calls for limiting the maximum size of the timestep one is ready to accept because, otherwise, the timesteps $\Delta_n t$ could reach unreasonable large values and then  the sampling of the solution could be too scarce. This phenomenology (and the consequent procedure of limiting the maximum allowed time interval) is parallel to that of the adaptive methods for ordinary differential equations \cite{Press2007}.

\begin{figure}
\begin{center}
\includegraphics[width=0.75\textwidth]{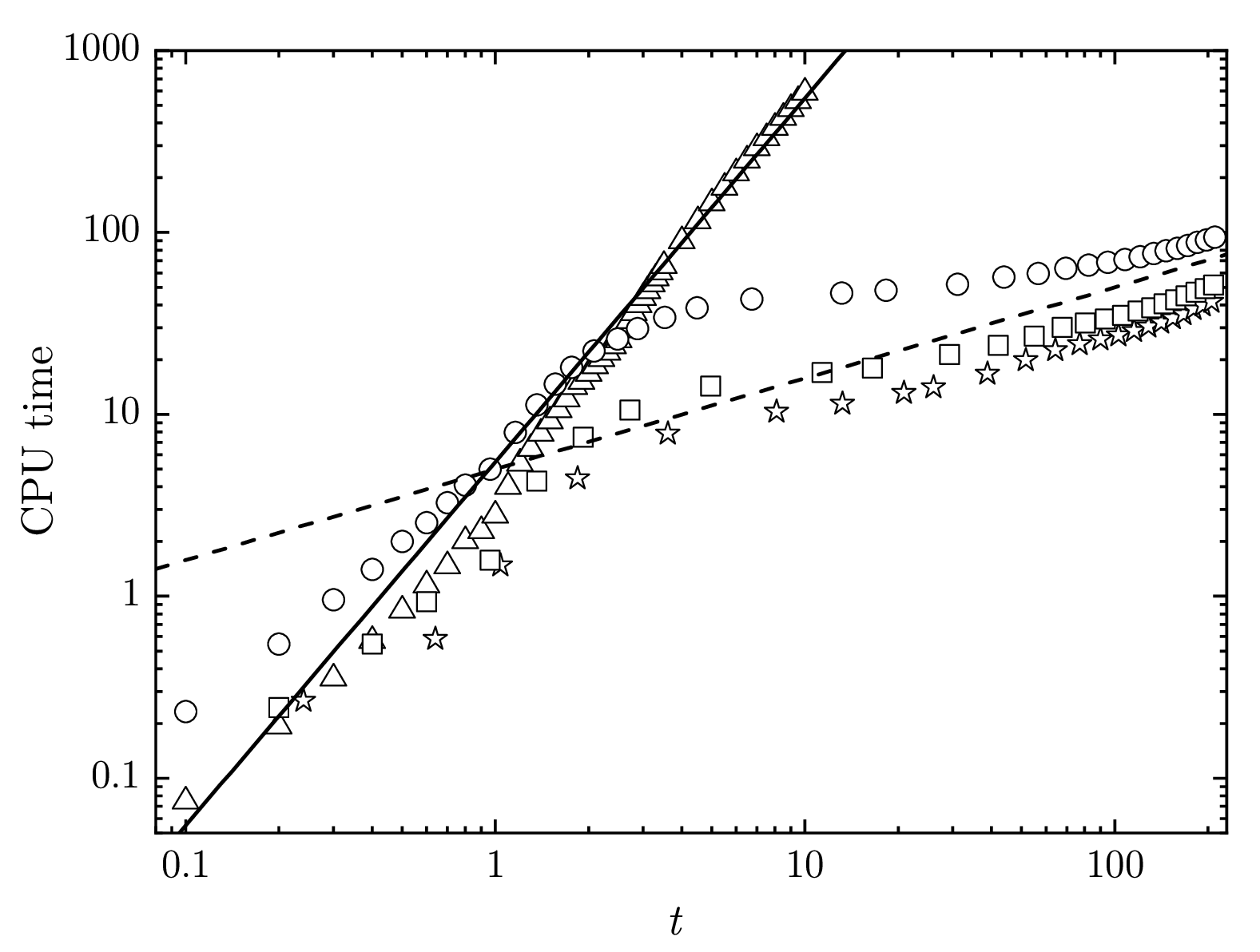}    
\caption{ Normalized computational time $T_\text{CPU}(t)$  required by the fixed-step method with $\Delta_n =0.01$ (triangles) and  by the adaptive method with $\tau_1=10^{-3}$ (stars), $\tau_2=5\times 10^{-4}$ (squares) and $\tau_3=10^{-5}$ (circles). They are averaged values over five runs. In all cases $\Delta x=\pi/40$. As reference, we have also plotted two lines corresponding to $T_\text{CPU}\sim t^2$ (solid line) and $T_\text{CPU}\sim t^{1/2}$ (dashed line). }
\label{fig:CPUtimesm2alpha1}
\end{center}
\end{figure}

It is also noticeable  that the number of symbols that appear in Fig.~\ref{fig:Ucase2} is larger for smaller values of the tolerance, which means that, in general, the larger  the tolerance, the larger  the timesteps.
This makes sense as one can use larger timesteps if one is ready to pay the price of larger numerical errors.  %
This can be clearly seen in the inset of  Fig.~\ref{fig:Ucase2} where the numerical error at the midpoint, $|u(x_j=1/2,t_n)-U_j^n|$, is plotted: typically, the larger the $\tau$ the larger the error.  For times larger than, say, $t=5$, the solution is close to its asymptotic value $u=2$ and the error becomes quite small. As reference, we have also plotted the numerical errors when a fixed timestep of size $\Delta_n=0.01$ is used. Except for small values of time, the errors of the adaptive method are similar (or smaller) than those of the method with fixed timestep, and, remarkably, this is reached with a huge improvement in the computation time.  This is clearly seen in Fig.~\ref{fig:CPUtimesm2alpha1} where we compare the CPU time required by the method with fixed timesteps and the method with adaptive timesteps (with three different values of the tolerance) to evaluate the numerical solutions. We see that for times $t\gtrsim 2$ the adaptive method is by far the best option.


The CPU times in Fig.~\ref{fig:CPUtimesm2alpha1}  are normalized times: they are not given in seconds but in units of $\mathcal{T}_{50}$, which is the CPU time employed by the method with fixed timesteps to get the solution of problem \eqref{eqsCaso1} when $50$ timesteps are used (in our computer $\mathcal{T}_{50} \approx 1.3$ seconds). Thus the normalized values $T_\text{CPU}(t)$ we report here  should be roughly independent of the computer one uses.   Notice that, as expected,  $T_\text{CPU} \propto t^2$ for the method with fixed timesteps (which comes from the fact that the computation time scale as the $N^2$, with  $N$ being  the number of timesteps). However, we see that the CPU time growths much more slower for the adaptive method: $T_\text{CPU}\sim t^{\beta}$  with  $\beta \approx 1/2$ for $t\gtrsim 2$.
This behaviour of  $T_\text{CPU}$ is very similar to the one observed for the case with constant $\gamma$  we studied in \cite{Yuste2016}, although the  relationship $T_\text{CPU}\sim t^{\beta}$ with $\beta$ constant is not as good here as we found  in \cite{Yuste2016} for the case with constant $\gamma$. This is a  consequence of the fact that  we consider here a variable order equation: in \cite{Yuste2016} we found that $\beta$ increases when so does $\gamma$, so that  a power fit for the CPU times with $\beta$ constant cannot be perfect if $\gamma$ varies with the time. Finally, we observe in Fig.~\ref{fig:CPUtimesm2alpha1} that, as expected, the adaptive method is faster when the tolerance increases.

\subsection{Case 2}

\begin{figure}
\begin{center}
\includegraphics[width=0.48\textwidth]{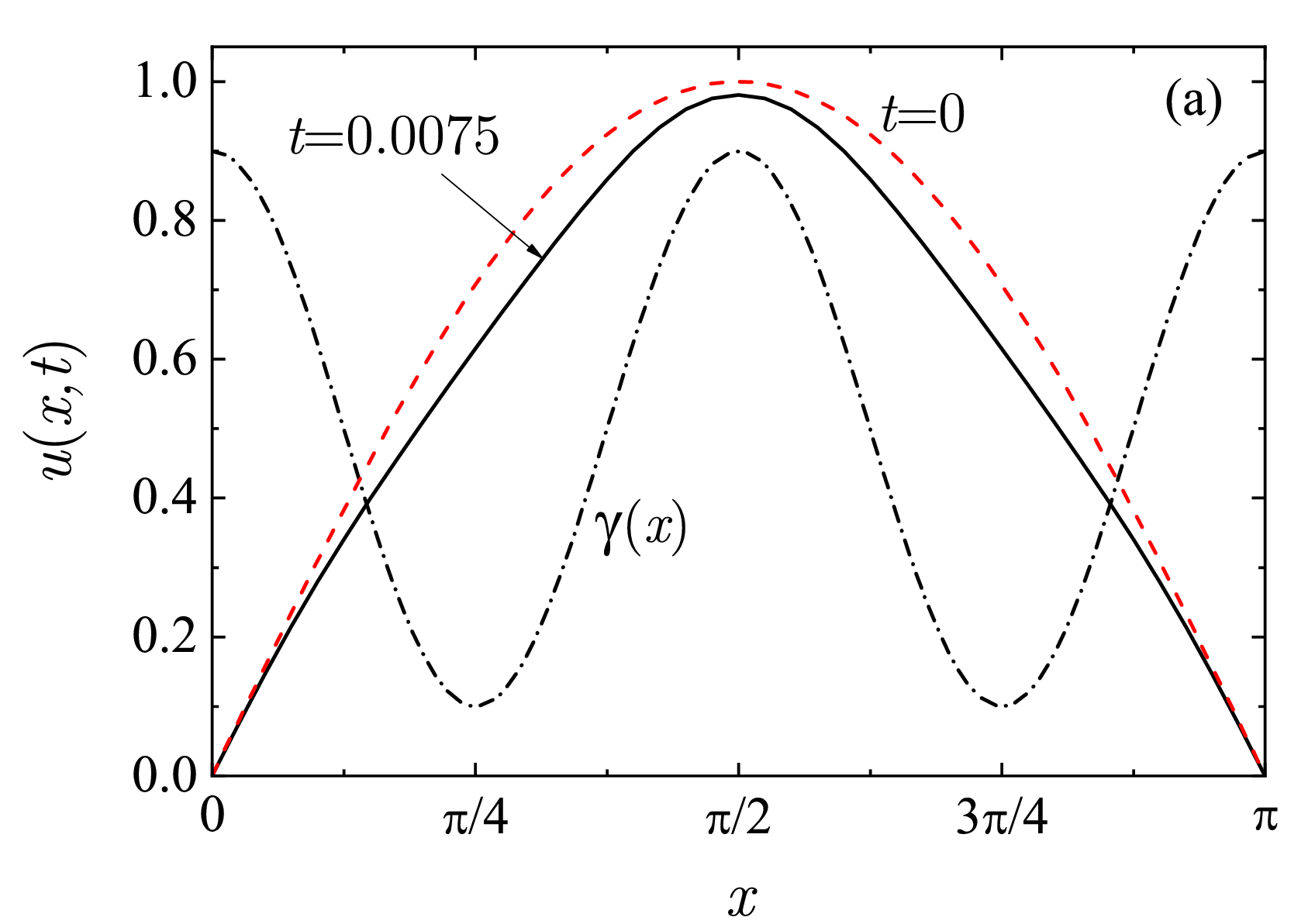} \hfill
\includegraphics[width=0.48\textwidth]{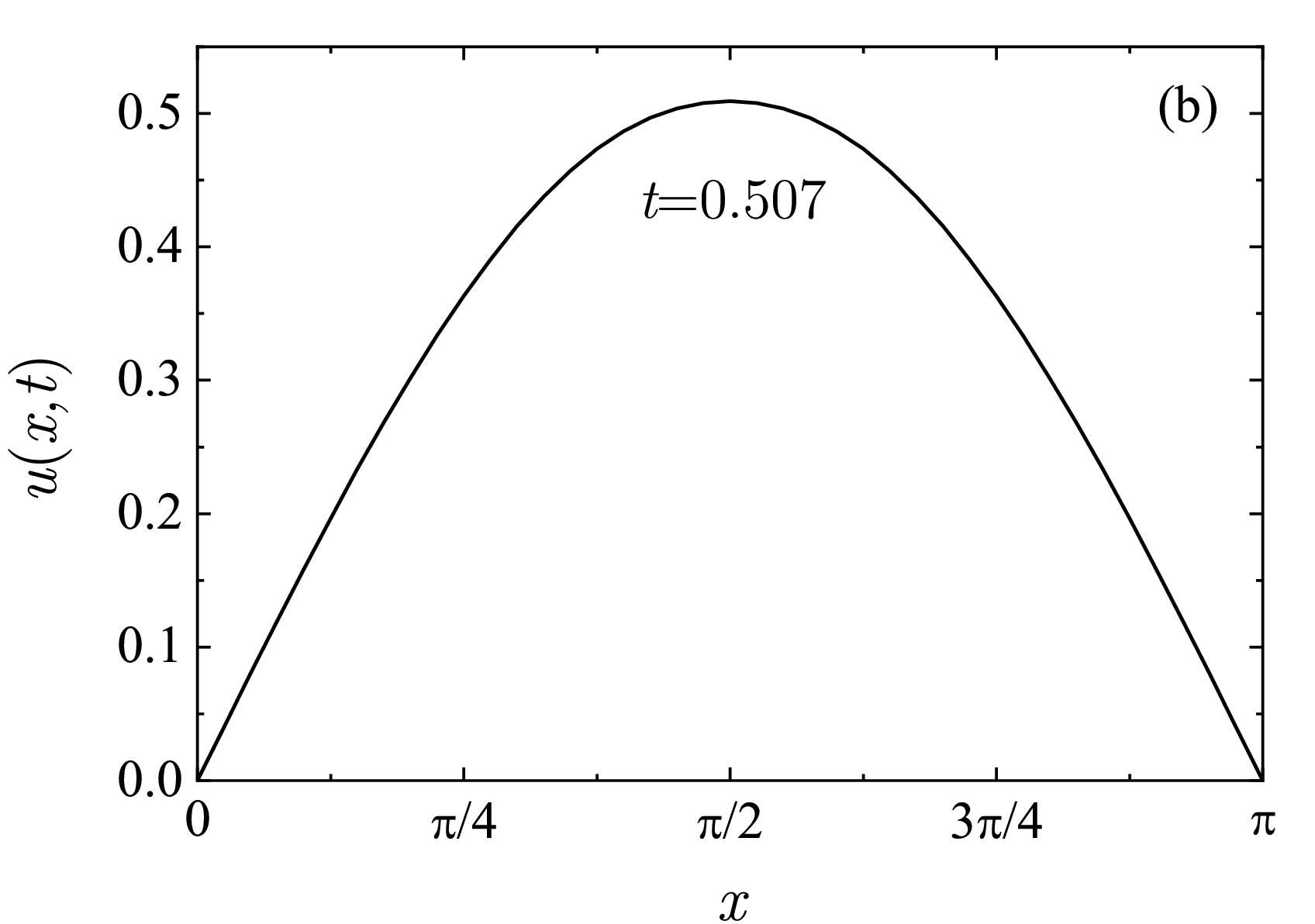}\\
\includegraphics[width=0.48\textwidth]{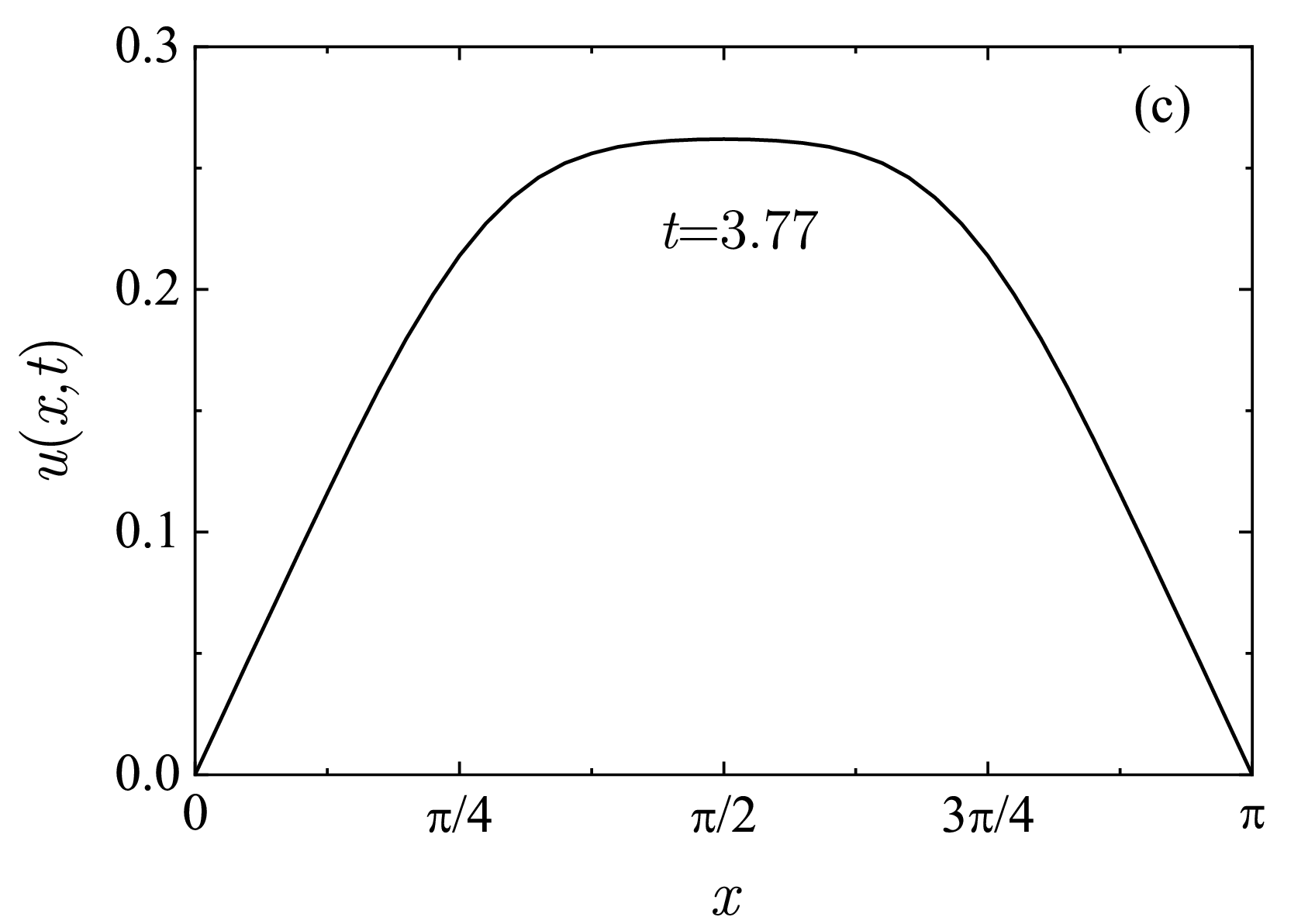} \hfill
\includegraphics[width=0.48\textwidth]{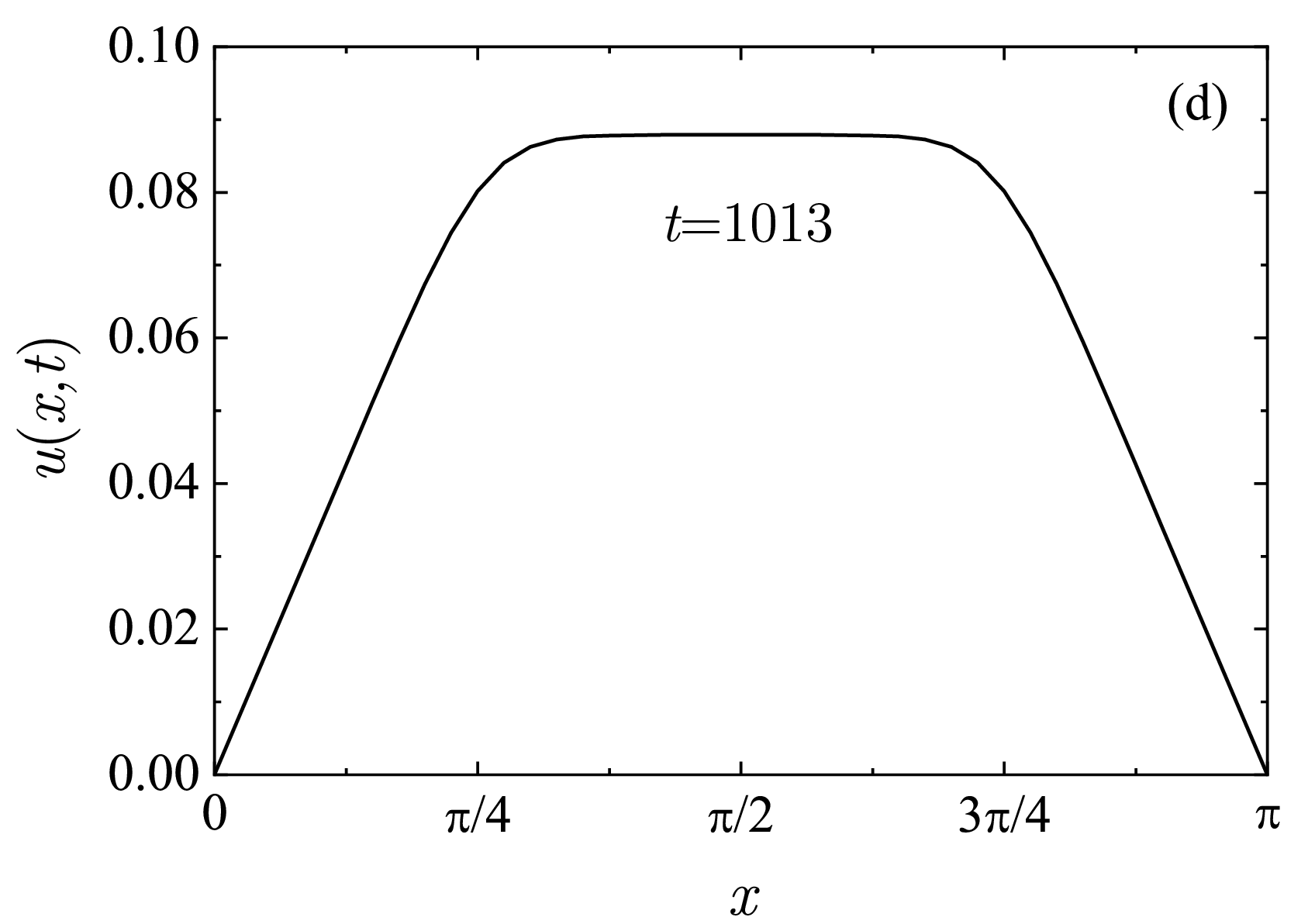}
\caption{Adaptive numerical solution $u(x, t)$ of the problem \eqref{eqsCaso2} (Case 2) with $\Delta x=\pi/40$ and $\tau=10^{-4}$ for times (a) $t=0.0075$, (b) $t=0.507$, (c) $t=3.77$ and (d) $t=1013$ (solid lines). As reference, in panel (a) we have also plotted $u(x,0)$ (dashed line) and the order $\gamma(x)$ of the fractional derivative  (dash-dotted line).  }
\label{fig:gammaCos2x}
\end{center}
\end{figure}

Now we consider a case where the order of the fractional derivative only depends on the position, i.e.,  $\gamma(x,t)=\gamma(x)$.
The problem is
\begin{subequations}
\label{eqsCaso2}
\begin{align}
        \frac{\partial^{\gamma}u(x,t)}{\partial t^\gamma} &=\frac{\partial^2 u(x,t)}{\partial x^2},  \label{eqCaso2} \\
        u(x,0)&=\sin x,\qquad 0\leq x \leq \pi, \label{eqCaso2CI} \\
        u(0,t)&=u(\pi,t)=0,  \label{eqCaso2CC}
\end{align}
\end{subequations}
where the derivative order is the periodic function  $\gamma(x)=\left[1+8\cos^2(2x)\right]/10$.  The solution obtained with the adaptive method with $\Delta x=\pi/40$ and tolerance $\tau=10^{-4}$ is shown in Fig.~\ref{fig:gammaCos2x} for four different times.  The solution is initially a sine function, $u(x,0)=\sin x$, but very quickly changes its form to a short of rounded triangle: see the panel for $t=7.5\times 10^{-3}$. After a while, the form of the solution is again similar to the initial one (see the panel for  $t=0.5075$). From here on the form of the solution goes towards curve similar to an isosceles trapezium with slowly decreasing height.
Although one has to go to times as large as $t\sim 10^3$  to see this with clarity (see the four panel), these large times can easily reached by our adaptive method; for example, only  134 timesteps were required to reach the time $t=1013$.
The fact that the central part of the solution becomes flatter and flatter as time increases is due to the dependence of $\gamma$ on $x$.
This example shows the usefulness of adaptive methods if the solution changes very quickly for very short times (so that one has to use very small timesteps to tracks down the solution in this time regime) and one has to go to long times (and then use large timesteps) to be able to appreciate the final form of the solution.  This example also shows that the VOFDE solution is qualitatively different from the corresponding FDE solution with gamma constant,  a case where the solution is always proportional to $\sin (x)$.

\begin{figure}
\begin{center}
\includegraphics[width=0.7\textwidth]{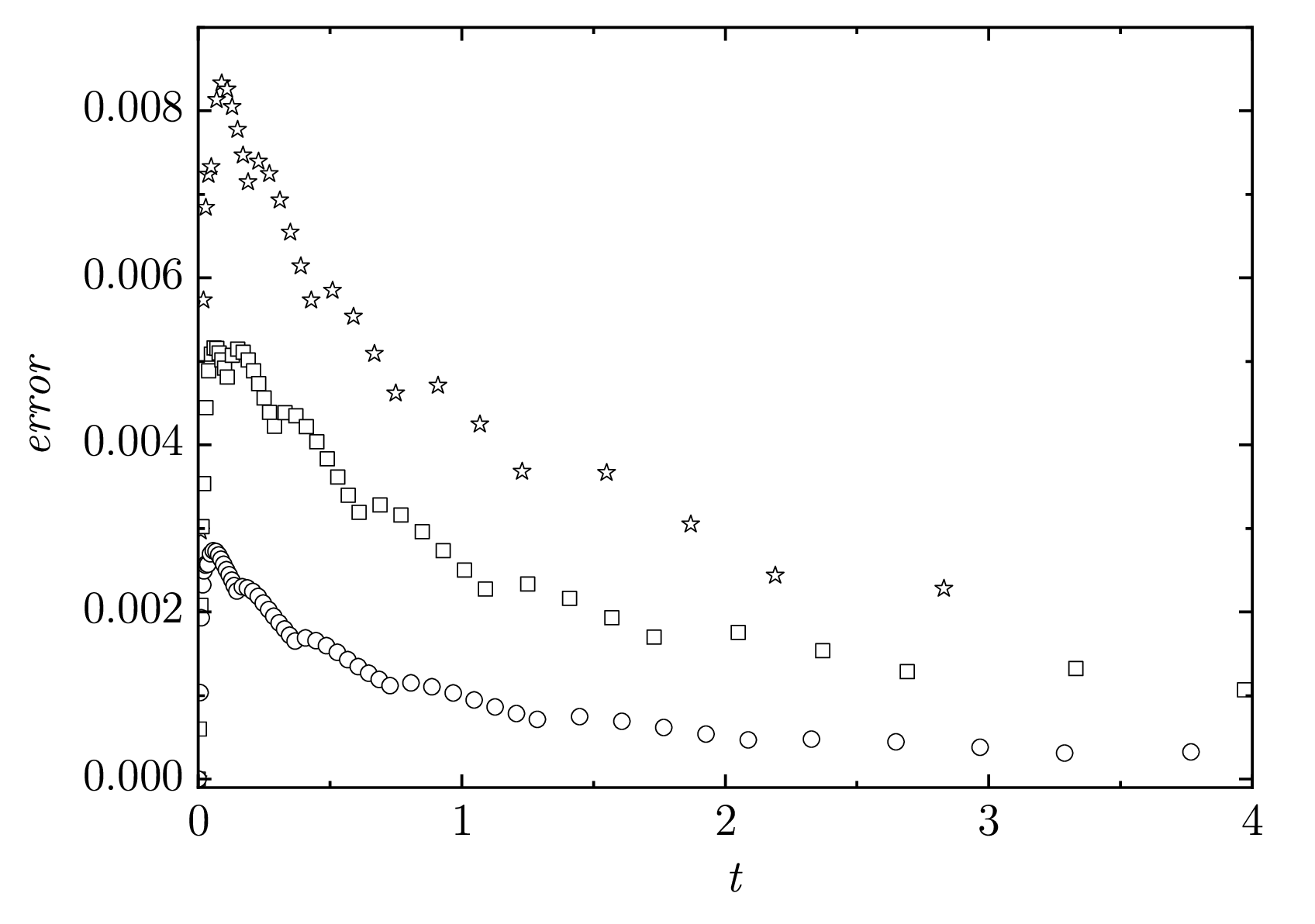}
\caption{
\sbytext{
Estimate  $|U_j^n-\widetilde{U}_j^n|$ of the error of the adaptive numerical solution $u(x, t)$ of the problem \eqref{eqsCaso2} (Case 2) when  $\Delta x=\pi/40$ and $\tau=10^{-3}$ (stars), $\tau=5\times 10^{-4}$ (squares), and $\tau=10^{-4}$ (circles).
}  }
\label{fig:ErrCase2}
\end{center}
\end{figure}

\sbytext{
Figure~\ref{fig:ErrCase2} shows the (estimated) error of our numerical method at the midpoint $x=\pi/2$. Since this VOFDE problem lacks an exact solution, we compute the error as $|U_j^n-\widetilde{U}_j^n|$, where $U_j^n$ is the estimated solution at the midpoint with $\Delta x=\pi/40$ and tolerances $\tau=10^{-3}$, $\tau=5\times 10^{-4}$, and $\tau=10^{-4}$. The quantity $\widetilde{U}_j^n$ corresponds to the numerical solution at the midpoint obtained with $\Delta x=\pi/160$ and $\tau=10^{-6}$. The rationale behind this procedure is that $\Delta x$ and $\tau$ are so small in the latter case that the true error $|U_j^n-u(x=\pi/2,t_n)|$ can be well approximated by $|U_j^n-\widetilde{U}_j^n|$.
We observe that the behavior of the error is very similar to that of Case 1 (see the inset in Fig.~\ref{fig:Ucase2}). In particular, we observe that the error is significantly smaller for smaller tolerances and that, after a short initial growth phase, it gradually decreases to values closer to the specified tolerance.
}

\subsection{Case 3}

Here we will also consider a case with a position dependent fractional order derivative  but now with nonhomogeneous boundary conditions, which implies that the solution goes, albeit very slowly, towards a nonzero stationary solution.

\begin{figure}
\begin{center}
\includegraphics[width=0.7\textwidth]{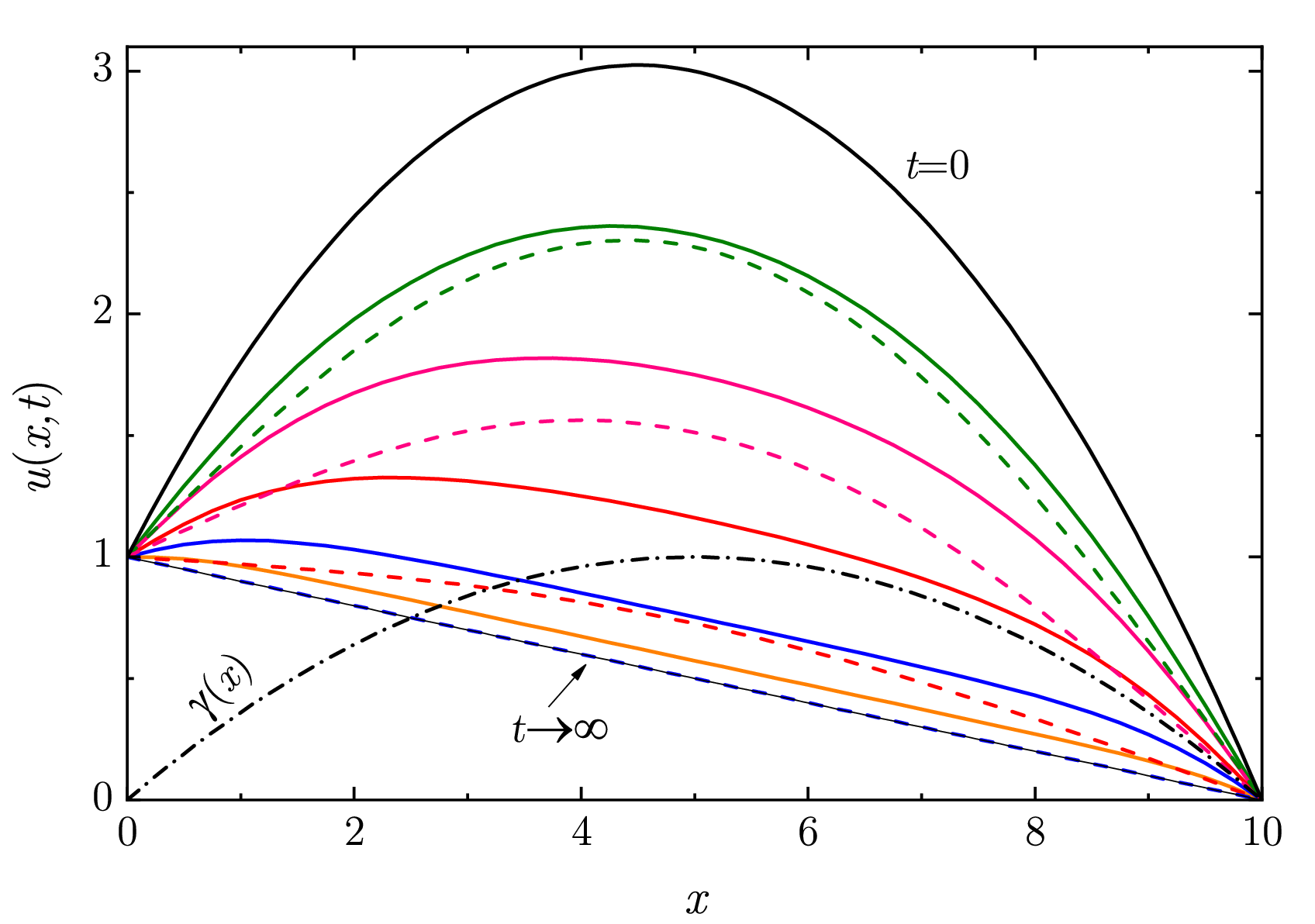}
\caption{Adaptive numerical solution $u(x, t)$ of the problem \eqref{eqsCaso3} (Case 3) with $\Delta x=0.01$ and $\tau=10^{-4}$. The solid lines represent the solution for, from top to bottom,  $t=0, 3.8, 9.6, 24.8, 93.0$ and $1350$ (black, green, pink, red, blue, and orange lines, respectively) for $\gamma(x)= (2/5)x(1-x/10)$ (dash-dotted line). The dashed lines are the solutions for the normal diffusion problem (i.e., $\gamma=1$) for, from top to bottom,  $t=3.8, 9.6, 24.9$ and $91.8$ (green, pink, red and blue lines, respectively). \sbytext{This last line overlaps the line corresponding to the final stationary solution, which for both $\gamma$ functions is the thin straight line going from $u=1$ to $u=0$.}}
\label{fig:gammaParabo}
\end{center}
\end{figure}

The problem we consider is
\begin{subequations}
\label{eqsCaso3}
\begin{align}
        \frac{\partial^{\gamma}u(x,t)}{\partial t^\gamma} &=\frac{\partial^2 u(x,t)}{\partial x^2}, \label{eqCaso3} \\
        u(x,0)&=(x+1)(1-x/10),\qquad 0\leq x \leq 10, \label{eqCaso3CI} \\
        u(0,t)&=1,\qquad u(10,t)=0, \label{eqCaso3CC}
\end{align}
\end{subequations}
with $\gamma=\gamma(x)=2x(1-x/10)/5$. The solution obtained with the adaptive method with $\Delta x=0.01$ and tolerance $\tau=10^{-4}$ is shown in Fig.~\ref{fig:gammaParabo} for several times.  We see that, initially, the solution changes very fast: the change of the solution from $t=0$ to $t\approx 3$ is similar to the change from $t\approx 25$ to $t\approx 90$.  This evolution is increasingly slower: note that even for $t=1350$ the solution has not yet reached the stationary solution. This behaviour is at odds with the behaviour of the solutions for the standard normal diffusion problem in which $\gamma=1$.  The solutions obtained by means of the  adaptive method for this case are also plotted in Fig.~\ref{fig:gammaParabo} (dashed lines). In particular, we see that the rates of convergence of the solutions towards the stationary solution  are vastly different: the solution for $\gamma=1$ overlaps the stationary solution already for $t\approx 90$  whereas the solution for $\gamma$ variable is clearly different from the stationary solution even for times as large as $t=1350$.   This phenomenology shows us once again the convenience of numerical methods with variable timesteps \sbytext{for solving} VOFDPEs.

\begin{figure}
\begin{center}
\includegraphics[width=0.7\textwidth]{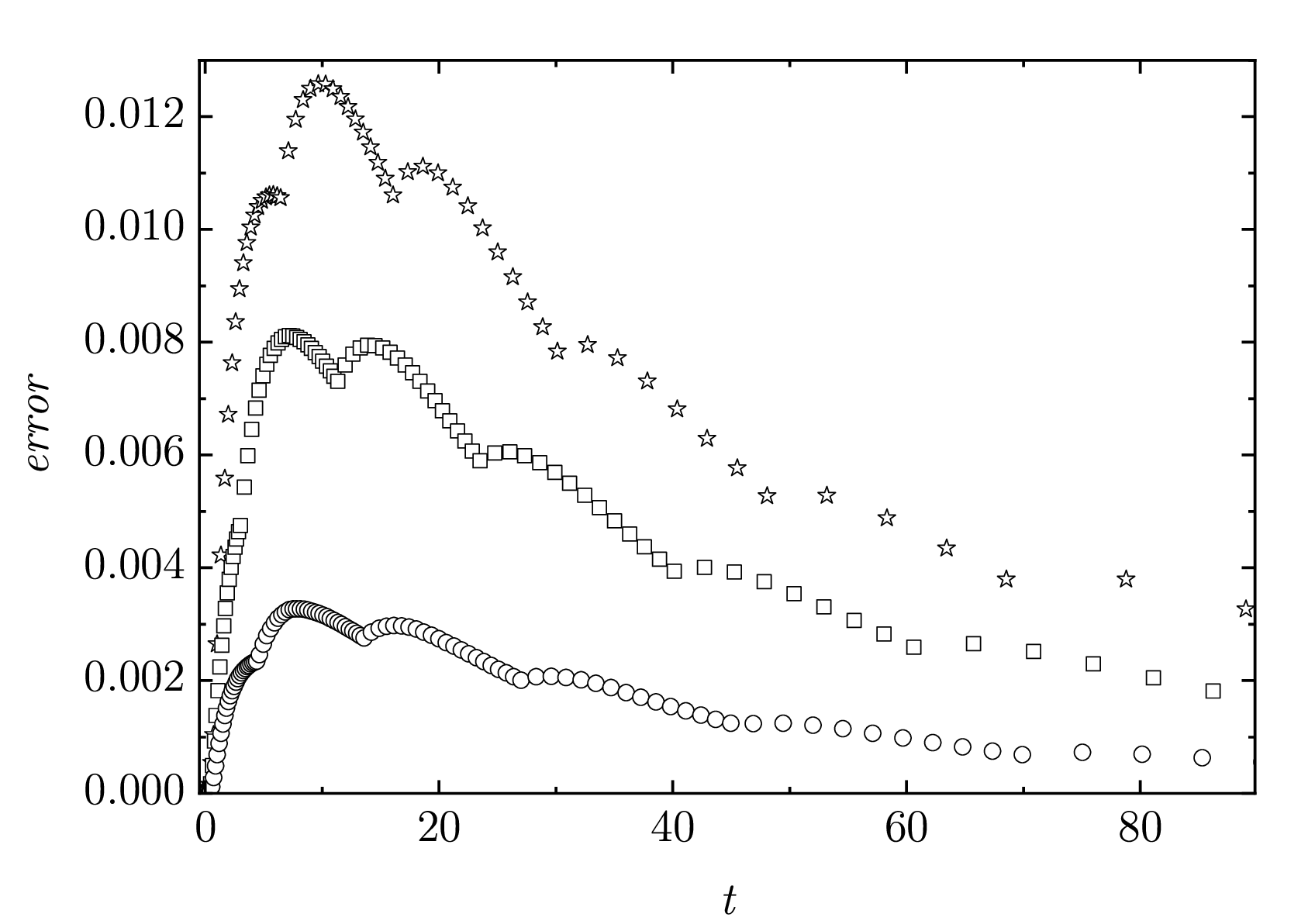}
\caption{
\sbytext{Estimate  $\max_{j}|U_j^n-\widetilde{U}_j^n|$ of the error of the adaptive numerical solution $u(x, t)$ of the problem \eqref{eqsCaso3} (Case 3) when  $\Delta x=\pi/40$ and $\tau=10^{-3}$ (stars), $\tau=5\times 10^{-4}$ (squares), and $\tau=10^{-4}$ (circles).
} }
\label{fig:ErrCase3}
\end{center}
\end{figure}

\sbytext{Figure~\ref{fig:ErrCase3} shows the (estimated) maximum error of our adaptive numerical method for several tolerances. Again, since this VOFDE problem lacks an exact solution, we estimate the error of the adaptive numerical algorithm as $\max_{j}|U_j^n-\widetilde{U}_j^n|$, where $U_j^n$ is the numerical solution for $\Delta x=1/4$ and tolerances $\tau=10^{-3}$, $\tau=5\times 10^{-4}$, and $\tau=10^{-4}$. The quantities $\widetilde{U}_j^n$ are the numerical solution when $\Delta x=1/16$ and $\tau=10^{-6}$. These latter values of $\Delta x$ and $\tau$ are so small that we can safely estimate the true error $|U_j^n-u(x=\pi/2,t_n)|$ by $|U_j^n-\widetilde{U}_j^n|$.
Again, we observe that the behavior of the error is very similar to that of Case 1 and Case 2. Errors are significantly smaller for small tolerances and have a tendency to decrease toward values closer to the tolerance.
}

\subsection{Case 4}

\begin{figure}
\begin{center}
\includegraphics[width=0.7\textwidth]{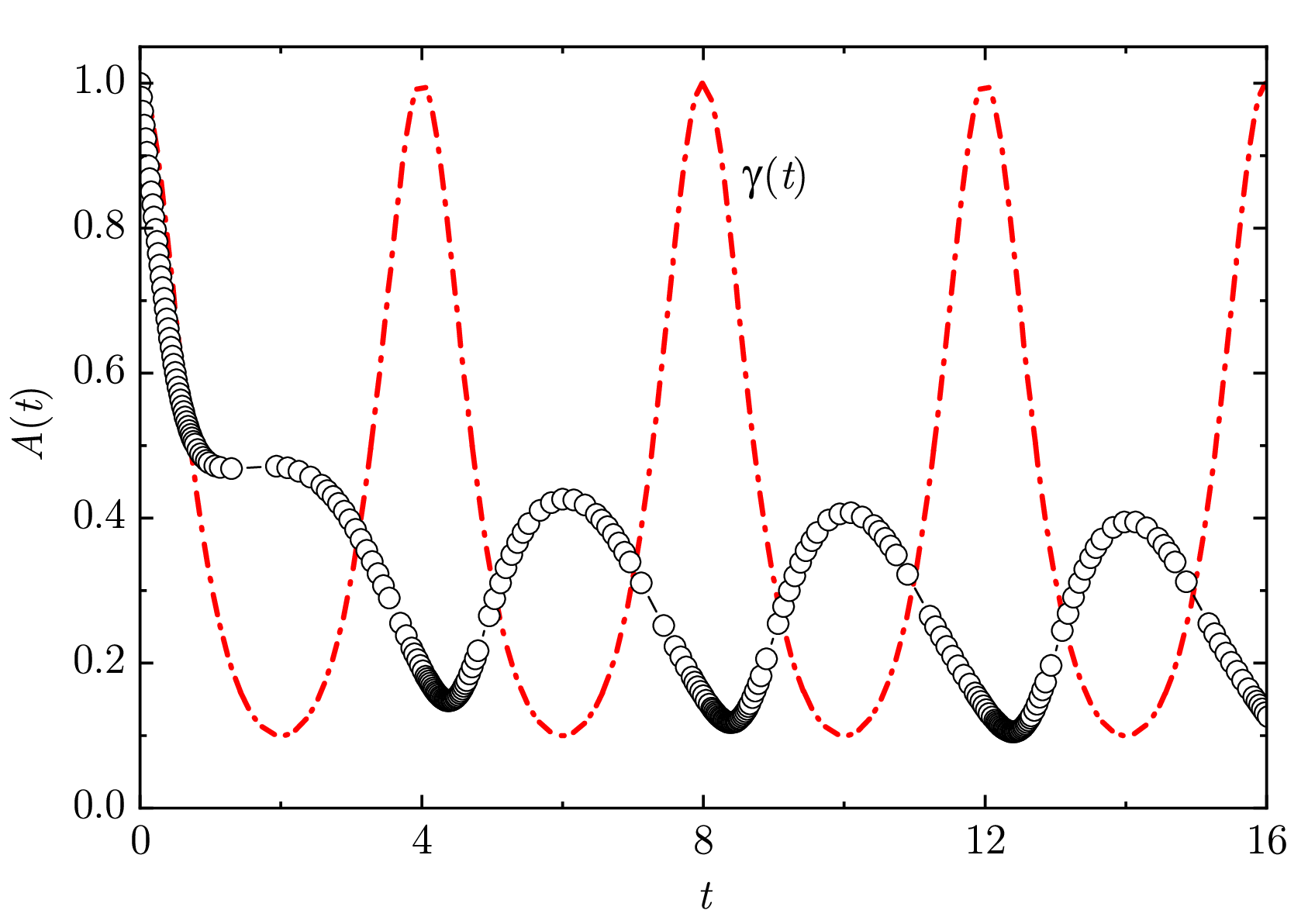}
\caption{Amplitude $A(t)$ of the solution $u(x, t)=A(t) \sin x$ of the problem \eqref{eqsCaso4} (Case 4). The points correspond to the numerical estimate $A(t_n)$ obtained by means of the present adaptive method with $\Delta x=\pi/40$ and $\tau=10^{-4}$. The continuous thin line through these points is an aid to the eye. For reference, we have also plotted the order of the fractional derivative $\gamma(t)$ (dash-dotted line).  }
\label{fig:gammaDenam}
\end{center}
\end{figure}

Now we consider a problem where the variable order fractional derivative $\gamma$ is a periodic function. The problem is
\begin{subequations}
\label{eqsCaso4}
\begin{align}
        \frac{\partial^{\gamma}u(x,t)}{\partial t^\gamma} &=\frac{\partial^2 u(x,t)}{\partial x^2},  \label{eqCaso4} \\
        u(x,0)&=\sin x,\qquad 0\leq x \leq \pi, \label{eqCaso4CI} \\
        u(0,t)&=1, \qquad u(\pi,t)=0, \label{eqCaso4CC}
\end{align}
\end{subequations}
where  $\gamma=\gamma(t)=\text{dn}(K(m) t/2,m)$, $m=0.99$, $K(\cdot)$ is the elliptical integral of first kind and $\text{dn}(\cdot)$ is the denam Jacobi elliptic function.  With this choice, $\gamma(t)$ is a periodic function of period 4 with flat valleys (i.e., with relatively small changes in the slope of the curve) around $t=2+4n$ with $n=0,1,\cdots$ and   sharp peaks around $4n$ with $n=0,1,\cdots$ (see Fig.~\ref{fig:gammaDenam}). The exact solution of this problem has the form $x(t)=A(t)\sin x$.  In  Fig.~\ref{fig:gammaDenam} we plot the numerical estimate of $A(t)$ obtained by the present adaptive method with $\Delta x=\pi/40$ and  tolerance $\tau=10^{-4}$. It is interesting to see how the amplitude of the solution syncs with $\gamma(t)$. This example shows us again that the adaptive method chooses the timesteps according to the behavior of the solution employing large timesteps when the curvature is small and short timesteps otherwise: notice that the points accumulate around $t\approx 4n$ while they scatter around $t\approx 2+4n$.  This implies an optimization of the number of timesteps required while maintaining the precision of the numerical solution. However, after the first few periods, the method is somewhat unsatisfactory as it continues to evaluate the optimal values of the timesteps period after period when they are essentially the same in each period. Of course, in these cases, one could stop this waste of computation time by stopping the adaptive procedure and just use (or at least use as the first trial timestep) the corresponding timestep from the previous period: $\Delta_n\approx\Delta_m$ with $t_n=t_m+T$, where $T$ is the period of the solution ($T=4$ in our example).

\begin{figure}
\begin{center}
\includegraphics[width=0.7\textwidth]{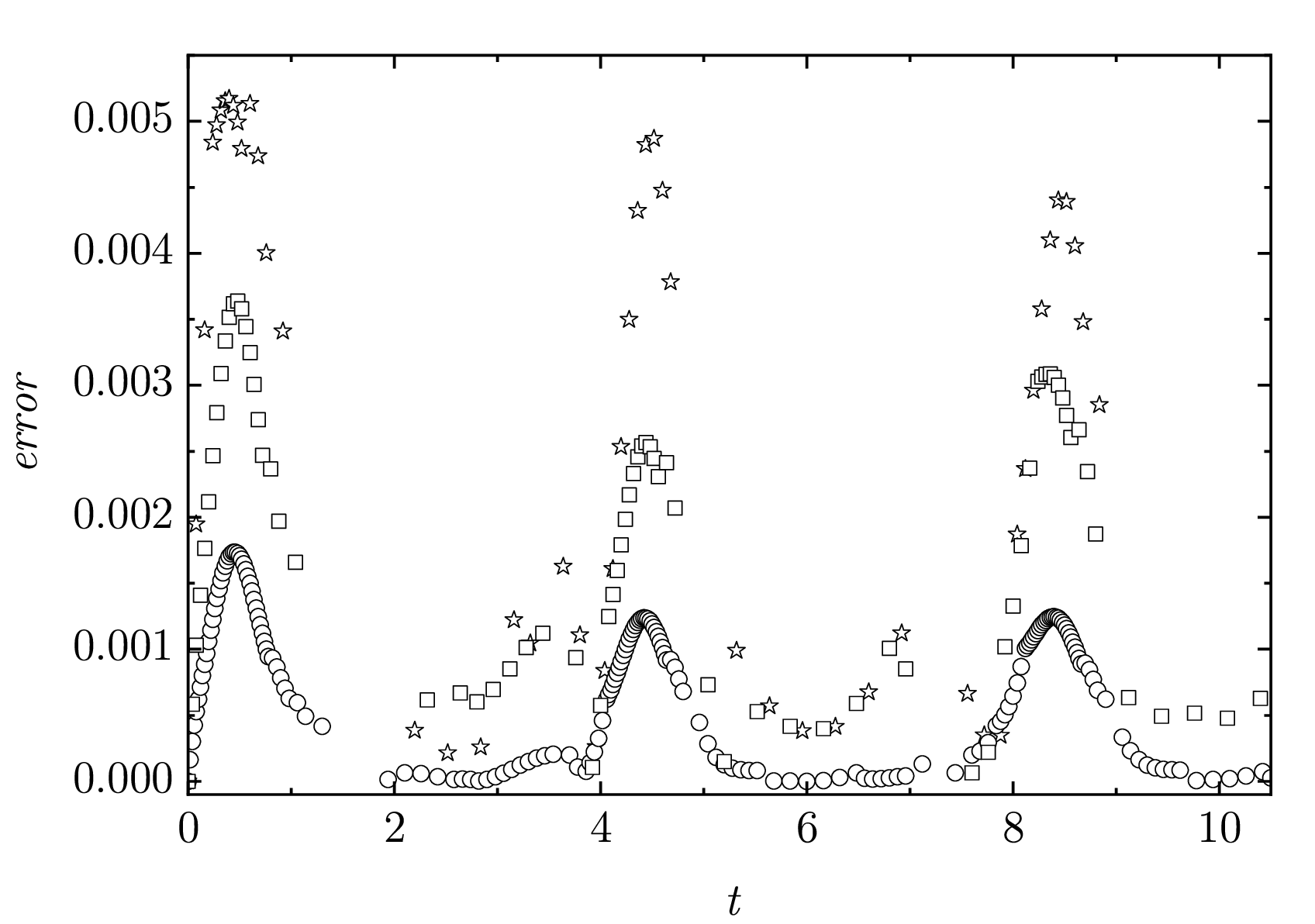}
\caption{
\sbytext{Estimate  $|U_j^n-\widetilde{U}_j^n|$ of the error at the midpoint  of the adaptive numerical solution $u(x, t)$ of the problem \eqref{eqsCaso4} (Case 4) when  $\Delta x=\pi/40$ and $\tau=10^{-3}$ (stars), $\tau=5\times 10^{-4}$ (squares), and $\tau=10^{-4}$ (circles).
}
 }
\label{fig:ErrCase4}
\end{center}
\end{figure}

\sbytext{
Figure~\ref{fig:ErrCase4} shows the (estimated) error of our adaptive numerical method at the midpoint $x=\pi/2$, i.e. the error in the value of $A(t)$. Since this VOFDE problem has no exact solution, we estimate the true error $|U_j^n-u(x=\pi/2,t_n)|$ by $|U_j^n-\widetilde{U}_j^n|$. Here $U_j^n$ is the numerical solution at the midpoint for $\Delta x=\pi/40$ and tolerances $\tau=10^{-3}$, $\tau=5\times 10^{-4}$, and $\tau=10^{-4}$, while $\widetilde{U}_j^n$ is the numerical solution at the midpoint for $\Delta x=\pi/160$ and $\tau=10^{-6}$.
Again, we see that the errors shrink significantly as the tolerance decreases, but now the time evolution of the error is different from that of cases 1--3. In those cases, the error went roughly monotonically to values close to the tolerances, but now the periodic nature of the solution is reflected in the roughly periodic behavior of the errors.
}

\section{Conclusion}
\label{secSumConclu}
A  previous finite-difference adaptive method \cite{Yuste2016} designed to solve some constant-order FDEs has been generalized to make it applicable to the variable-order version of these FDEs [see Eq.~\eqref{difuEqCaputo}].  The method has two main ingredients.  First, a finite difference scheme able to work with variable timesteps. In this paper we have used a scheme based on the L1 discretization of the Caputo fractional time derivative as well as the standard three-point centered discretization formula of the spatial derivative. The scheme is unconditionally stable. The second ingredient is a algorithm for choosing the size of the timesteps.  Here we have used a step-doubling algorithm to keep the numerical truncation error to values of the order of a quantity (the tolerance) that we preset.
A  characteristic of the FDEs is that their solutions typically show very disparate rates of change, generally with very fast changes for small times and very slow changes (``aging'') for large times.
\sbytext{The presence of very different temporal regimes is particularly pronounced in the FDEs with variable order $\gamma$. This is because, in addition to the already characteristic behavior of the rates of change of the FDEs mentioned above, there is the temporal and spatial variability resulting from the temporal and spatial dependence of the order $\gamma(x,t)$ of the FDEs.
}
For this type of behavior \sbytext{with such different time scales}, finite difference methods with fixed timesteps are doomed to be very inefficient: either they overlook what the short-time solution is like if, to describe the long-time solutions, a very large the timestep is chosen, or they cannot describe the solution for long times if, to describe the solutions for short times, the timestep is chosen too small.
The adaptive difference method we have presented overcomes these difficulties by changing the size of the timesteps according to the behavior (the rate of change) of the solutions.   In this way the solution of the problem is not only obtained with an accuracy that can be chosen freely, but is also achieved very efficiently with CPU times that, in many cases, are much smaller than the CPU times required by the corresponding standard method that uses a constant timestep.

The present method can be extended in at least two directions related to its two main ingredients. Another finite difference scheme with variable timesteps could be used by replacing the L1 discretization by another discretization formula (e.g., by the L2, L2-1$_\sigma$ or L1-2 formulas \cite{Alikhanov2015,Gao2014,Oldham1974}) and, also,  a different algorithm could be used to choose the size of the timesteps  \cite{Eriksson2017,Jannelli2020,Press2007}.
Work is in progress along these lines.


\subsection*{Acknowledgments}
S.B.Yuste acknowledges financial support
\sbytext{from Grant PID2020-112936GB-I00 funded by MCIN/AEI/10.13039/501100011033, and from Grant IB20079 funded by Junta de Extremadura (Spain) and by ERDF A way of making Europe.}



\end{document}